\documentstyle[11pt]{article}
\newcommand{\de}{\delta}
\newcommand{\G}{\Gamma}
\newcommand{\beq}[1]{\begin{equation}\label{#1}}
\newcommand{\eq}{\end{equation}}
\newcommand{\beqn}[1]{\begin{eqnarray}\label{#1}}
\newcommand{\eqn}{\end{eqnarray}}
\newcommand{\p}{\partial}
\newcommand{\mbesj}{I_\nu^{(j)}((1-q^2)z;q^2)}
\newcommand{\mbes}{I_\nu^{(1)}((1-q^2)z;q^2)}

\newcommand{\mmbesj}{I_{-\nu}^{(j)}((1-q^2)z;q^2)}
\newcommand{\mbess}{I_\nu^{(2)}((1-q^2)z;q^2)}
\newcommand{\mbesse}{I_\nu^{(3)}((1-q^2)z;q^2)}
\newcommand{\mmbesse}{I_{-\nu}^{(3)}((1-q^2)z;q^2)}
\newcommand{\makd}{K_\nu^{(1)}((1-q^2)z;q^2)}
\newcommand{\makdo}{K_\nu^{(2)}((1-q^2)z;q^2)}
\newcommand{\makdon}{K_\nu^{(3)}((1-q^2)z;q^2)}
\newcommand{\makdj}{K_\nu^{(j)}((1-q^2)z;q^2)}
\newcommand{\qqexp}{e_{q^2}(\frac{(1-q^2)^2}{4}z^2)}
\newcommand{\qqeexp}{E_{q^2}(-\frac{(1-q^2)^2}{4}q^2z^2)}
\newtheorem{predl}{Proposition}[section]
\newtheorem{defi}{Definition}[section]
\newtheorem{rem}{Remark}[section]
\newtheorem{cor}{Corollary}[section]

\begin{document}

\vspace{10mm}
\begin{flushright}
 ITEP-TH-56/00\\
\end{flushright}
\vspace{10mm}
\begin{center}
{\Large \bf $q$-Bessel-Macdonald functions}\\
\vspace{5mm}
V.-B.K.Rogov \footnote{The work was supported by the Russian Foundation for
Fundamental Research (grant no. 0001-00143) and the NIOKR MPS RF}\\
101475, Moscow, MIIT\\
e-mail vrogov@cemi.rssi.ru
\end{center}

\begin{abstract}

The modified $q$-Bessel functions and the $q$-Bessel-Macdonald functions
of the first and second kind are introduced. Their definition is based on
representations as power series. Recurrence
relations, the $q$-Wronskians, asymptotic decompositions and
$q$-integral representations are received. In addition, the $q$-Bessel-Macdonald
function of kind 3 is determined by its $q$-integral representation.
\end{abstract}

\section{Introduction}
\setcounter{equation}{0}

The $q$-Bessel-Jackson functions of kinds 1, 2 and 3 were introduced at
the very beginning of the century \cite{Ja}. Their properties were considered in
\cite{KS,Is,V,VK}.  The definitions of the modified $q$-Bessel-Jackson
functions ($q$-MBF) and $q$-Macdonald functions ($q$-MF) and their
properties were given in \cite{OR1}.

The $q$-analogs of the modified Bessel functions and Macdonald functions
are interesting because they arise in the harmonic analysis on the
quantum homogeneous spaces as in the classical case. If we suppose that the
commutation relations in the universal enveloping algebra and the
commutation relations for the generators of the quantum Lobachevsky space
\cite{OR0}
are determined by different parameters than the eigenvalue problem for the
second Casimir operator leads us to the difference equation depending on
one parameter. This equation is the second order one for three values of
this parameter. These values of parameter correspond the
$q$-Bessel-Jackson functions of kind 1 $J_\nu^{(1)}$, kind 2 $J_\nu^{(2)}$,
and kind 3 $J_\nu^{(3)}$. This fact allows to consider these functions with
uniform point of view. Same authors name the $q$-Bessel-Jackson functions
of kind 3 by the Hahn-Exton ones

In this work we start from the definition of the $q$-MBF as a solution of
the second order difference equation. The $q$-MBFs are connected with
$q$-Bessel functions as in the classical case. We determine the actions of
the difference operators, the recurrence relations, and the $q$-Wronskians
for them. The Laurent series for the $q$-MBF of kinds 1 and 2 are very
important for a determination of the $q$-MFs. At last we
represent the $q$-MFs of kinds 1 and 2 by the Jackson
$q$-integral. We determine the $q$-MFs of kind 3 by their $q$-integral
representation, and then we receive the expression of these functions by
the $q$-MBFs of kind 3.
\newpage

\section{Preliminary results}
\setcounter{equation}{0}

{\bf 2.1. The modified $q$-Bessel functions}

\vspace{5mm}

In \cite{Ja} the $q$-Bessel functions were defined as follows:
\beq{2.1}
J_\nu^{(1)}(z,q)=\frac{(q^{\nu+1},q)_\infty}{(q,q)_\infty}(z/2)^\nu
\phantom1_2\Phi_1(0,0;q^{\nu+1};q,-\frac{z^2}{4}),
\eq
\beq{2.2}
J_\nu^{(2)}(z,q)=\frac{(q^{\nu+1},q)_\infty}{(q,q)_\infty}(z/2)^\nu
\phantom1_0\Phi_1(-;q^{\nu+1};q,-\frac{z^2q^{\nu+1}}{4}),
\eq
\beq{2.3}
J_\nu^{(3)}(z,q)=\frac{(q^{\nu+1},q)_\infty}{(q,q)_\infty}(z/2)^\nu
\phantom1_1\Phi_1(-;q^{\nu+1};q,-\frac{z^2q^{\frac{\nu+1}2}}{4}).
\eq
where $\phantom1_r\Phi_s$ is basic hypergeometric function \cite {GR},
$$
\phantom1_r\Phi_s(a_1,\cdots,a_r;b_1,\cdots,b_s;q,z)=
\sum_{n=0}^\infty\frac{(a_1,q)_n\ldots(a_r,q)_n}
{(q,q)_n(b_1,q)_n\ldots(b_s,q)_n}[(-1)^nq^{n(n-1)/2}]^{1+s-r}z^n.
$$

It allows to introduce the modified $q$-Bessel functions
($q$-MBFs) using (\ref{2.1}), (\ref{2.2}) and (\ref{2.3}) similarly to
the classical case \cite{BE}.

\begin{defi}\label{d2.1}
The modified $q$-Bessel functions are the functions
$$
I_\nu^{(j)}(z,q)=\frac{(q^{\nu+1},q)_\infty}{(q,q)_\infty}(z/2)^\nu
\phantom1_\de\Phi_1\left(\underbrace{0,\ldots,0}_{\mbox{$\de$}};q^{\nu+1};
q,\frac{z^2q^{\frac{\nu+1}2(2-\de)}}{4}\right).
$$

Here
\beq{2.4}
\de=\left\{
\begin{array}{lcl}
2 & {\rm for} & j=1\\
0 & {\rm for} & j=2\\
1 & {\rm for} & j=3.\\
\end{array}
\right.
\eq
\end{defi}
Obviously,
$$
I_\nu^{(j)}(z,q)=e^{-\frac{i\nu\pi}{2}}J_\nu^{(j)}(e^{i\pi/2}z,q),
\qquad j=1,2,3.
$$
In the sequel we consider the functions
\beq{2.5}
I_\nu^{(1)}((1-q^2)z;q^2)=\sum_{k=0}^\infty\frac{(1-q^2)^k(z/2)^{\nu+2k}}
{(q^2,q^2)_k\Gamma_{q^2}(\nu+k+1)},\qquad |z|<\frac{2}{1-q^2},
\eq
\beq{2.6}
I_\nu^{(2)}((1-q^2)z;q^2)=
\sum_{k=0}^\infty\frac{q^{2k(\nu+k)}(1-q^2)^k(z/2)^{\nu+2k}}
{(q^2,q^2)_k\Gamma_{q^2}(\nu+k+1)}.
\end{equation}
\beq{2.7}
I_\nu^{(3)}((1-q^2)z;q^2)=
\sum_{k=0}^\infty\frac{q^{k(\nu+k)}(1-q^2)^k(z/2)^{\nu+2k}}
{(q^2,q^2)_k\Gamma_{q^2}(\nu+k+1)}.
\eq
If $|q|<1,$ the series (\ref{2.5}) and (\ref{2.6}) are absolutely
convergent for all $z\ne0$. Consequently, $\mbess$ and
$I_\nu^{(3)}((1-q^2)z;q^2)$ are holomorphic functions outside a
neighborhood of zero.
\begin{rem}\label{r2.1}
$$
\lim_{q\to1-0}I_\nu^{(j)}((1-q^2)z;q^2)=I_\nu(z),\qquad j=1, 2, 3.
$$
\end{rem}
\begin{predl}\label{p2.1}
The function $\mbesj$ is a solution of the difference equation
\beq{2.8}
f(q^{-1}z)-(q^{-\nu}+q^\nu)f(z)+f(qz)=
q^{-\de}\frac{(1-q^2)^2}4z^2f(q^{1-\de}z).
\eq
$j=1, 2, 3$ are connected with $\de=2, 0, 1$ by relations (\ref{2.4})
\end{predl}
\begin{cor}\label{c2.1}
The function $I_{-\nu}^{(j)}((1-q^2)z;q^2)$ satisfies equation
(\ref{2.8}).
\end{cor}
\begin{predl}\label{p2.2}
The functions $I_\nu^{(j)}((1-q^2)z;q^2)$ satisfies the relations
$$
\frac2{(1+q)z}\p_q{z^\nu I_{-\nu}^{(j)}((1-q^2)z;q^2)}=
q^{-\frac{2-\delta}2(\nu-1)}z^{\nu-1}
I_{-\nu+1}^{(j)}((1-q^2)q^{\frac{2-\delta}2}z;q^2),
$$
$$
\frac{2}{(1+q)z}\p_qz^\nu I_\nu^{(j)}((1-q^2)z;q^2)=
q^{-\frac{2-\delta}2(\nu-1)}z^{\nu-1}
I_{\nu-1}^{(j)}((1-q^2)q^{\frac{2-\delta}2}z;q^2),
$$
where the operator $\p_q$ is defined as $\p_qf(z)=\frac{f(z)-f(qz)}
{(1-q)z}$.
\end{predl}
\begin{predl}\label{p2.3}
The functions $I_\nu^{(j)}((1-q^2)z;q^2)$ satisfy the recurrence relations
$$
q^{-\frac{2-\delta}2\nu}I_{\nu-1}^{(j)}((1-q^2)z;q^2)-
q^{\frac{2-\delta}2\nu}I_{\nu+1}^{(j)}((1-q^2)z;q^2)=
$$
$$
=\frac{2}{(1-q^2)z}(q^{-\nu}-q^\nu)I_\nu^{(j)}((1-q^2)q^{\frac\delta2}z;q^2),
$$
$$
q^{-\frac{2-\delta}2\nu}I_{\nu-1}^{(j)}((1-q^2)z;q^2)+
q^{\frac{2-\delta}2\nu}I_{\nu+1}^{(j)}((1-q^2)z;q^2)
=\frac{4}{(1-q^2)z}I_\nu^{(j)}((1-q^2)q^{-\frac{2-\delta}2}z;q^2)-
$$
$$
-\frac{2}{(1-q^2)z}(q^{-\nu}+q^\nu)I_\nu^{(j)}((1-q^2)q^{\frac\delta2}z;q^2).
$$
\end{predl}

\begin{defi}\label{d2.2}
The $q$-Wronskian of two solutions $f_\nu^1(z)$ and $f_\nu^2(z)$ of a
second-order difference equation is defined as follows:
$$
W(f_\nu^1,f_\nu^2)(z)=f_\nu^1(z)f_\nu^2(qz)-f_\nu^1(qz)f_\nu^2(z).
$$
\end{defi}

If the $q$-Wronskian does not vanish, then any solution of the
second-order difference equation can de written in form
$$
f_\nu(z)=C_1f_\nu^1(z)+C_2f_\nu^2(z).
$$
In this case the functions $f_\nu^1(z)$ and $f_\nu^2(z)$ form a fundamental
system of the solutions of the given equation.

\begin{predl}\label{p2.4}
If $\nu\ne k+\frac{i\pi m}{\ln q}, ~k$ and $m$ are integers, then the
functions $\mbesj$ and $I_{-\nu}^{(j)}((1-q^2)z;q^2)$ form a fundamental
system of the solutions of equation (\ref{2.8})
($z\ne\pm\frac{2q^{-r}}{1-q^2}, r=0, 1,\ldots$).
\end{predl}

This Proposition is following from
\beq{2.9}
W(z)=\left\{
\begin{array}{lcl}
\frac{q^{-\nu}(1-q^2)}{\Gamma_{q^2}(\nu)\Gamma_{q^2}(1-\nu)}
\qqexp & {\rm for}&\de=2\\
\frac{q^{-\nu}(1-q^2)}{\Gamma_{q^2}(\nu)\Gamma_{q^2}(1-\nu)}&
{\rm for}&\de=1\\
\frac{q^{-\nu}(1-q^2)}{\Gamma_{q^2}(\nu)\Gamma_{q^2}(1-\nu)}
\qqeexp &{\rm for}&\de=0.\\
\end{array}
\right.
\eq
Obviously, this function is defined for $z\ne\pm\frac{2q^{-r}}{1-q^2},
r=0, 1,\ldots $ and does not vanish.

If $\nu=n$ is an integer, then in view of (\ref{2.5}) - (\ref{2.7})
$$
I_{-n}^{(j)}((1-q^2)z;q^2)=I_n^{(j)}((1-q^2)z;q^2),\quad j=1, 2, 3.
$$

The following relations take place
$$
\mbes=\qqexp\mbess.
$$
$$
\mbess=\qqeexp\mbes.
$$
\begin{predl}\label{p2.5}
The function $\mbes$ is a meromorphic function outside a\\ neighborhood of
zero, with simple poles at the points $z=\pm\frac{2q^{-r}}{1-q^2}, r=0,
1,\ldots$
\end{predl}
\begin{rem}\label{r2.2}
If $q\to1-0$, then the poles of the function $\mbes$
$$
z_r=\pm\frac{2q^{-r}}{1-q^2},\qquad r=0, 1,\ldots
$$
tend to infinity along the real axis.
\end{rem}

{\bf 2.2. Laurent type series for $q$-MBFs of kinds 1 and 2}

\begin{predl}\label{p2.6}
For $z\ne0$, and $j=1, 2$ ~$q$-MBF can be represented by form
$$
I_\nu^{(j)}((1-q^2)z;q^2)=\frac{a_\nu}{\sqrt{z}}
\left[\phi^{(j)}(z)\Phi_\nu(z)+
ie^{i\nu\pi}\phi^{(j)}(-z)\Phi_\nu(-z)\right],
$$
where
$$
\phi^{(j)}(z)=\left\{
\begin{array}{lcl}
e_q\left(\frac{1-q^2}2z\right) & {\rm for} & j=1\\
E_q\left(\frac{1-q^2}2z\right) & {\rm for} & j=2,\\
\end{array}
\right.
~~\Phi_\nu(z)=\phantom._2\Phi_1
\left(q^{\nu+\frac12},q^{-\nu+\frac12};-q;q,\frac{2q}{(1-q^2)z}\right),
$$
\beq{2.10}
a_\nu=\sqrt{\frac2{1-q^2}}\frac{I_\nu^{(2)}(2;q^2)}
{(-1,q)_\infty\Phi_\nu(\frac2{1-q^2})}.
\eq
\end{predl}
\begin{predl}\label{p2.7}
The coefficients $a_\nu$ satisfy the relations
\beq{2.11}
a_{\nu+1}=a_\nu q^{-\nu-\frac12},
\eq
\beq{2.12}
a_\nu a_{-\nu}=\frac{q^{-\nu+\frac12}}
{2\G_{q^2}(\nu)\G_{q^2}(1-\nu)\sin\nu\pi},
\eq
\end{predl}
If $\nu=n$, then from (\ref{2.11}) we have
$$
a_n=a_{n-k}q^{-k/2(2n-k)(\frac32-\frac{|\de-1|}2)}.
$$
Assume $k=2n$. Then $a_n=a_{-n}$. From (\ref{2.12}) and
\cite{GR} (1.10.6) we get
$$
(a_n)^2=\frac{q^{-n^2+1/2}\ln{q^{-2}}}{2\pi(1-q^2)}.
$$

{\bf 2.3. The $q$-Macdonald functions of kinds 1 and 2}

\begin{defi}\label{d2.3}
We define the $q$-Macdonald function ($q$-MF) for $j=1, 2, 3$ as follows:
$$
\makdj=
$$
\beq{2.13}
=\frac12q^{-\nu^2+\nu}\G_{q^2}(\nu)\G^{q^2}(1-\nu)
\left[\sqrt{\frac{a_\nu}{a_{-\nu}}}\mmbesj-
\sqrt{\frac{a_{-\nu}}{a_\nu}}\mbesj\right],
\eq
where $a_\nu^{(j)}$ is determined by (\ref{2.10}).
\end{defi}

As in the classical case, this definition must be extended to integral
values of $\nu=n$ by passing to the limit in (\ref{2.13}).

\begin{predl}\label{p2.8}
The $q$-MF $\makd$ is represented by form
$$
\makd=\frac{q^{-\nu^2+1/2}}{2\sqrt{a_\nu a_{-\nu}}\sqrt{z}}
e_q\left(-\frac{1-q^2}2z\right)\Phi_\nu(-z),
$$
and hence it is a holomorphic function in the region
$\Re e z>\frac{2q}{1-q^2}$.
\end{predl}

\begin{predl}\label{p2.9}
The $q$-MF $\makdo$ is represented by form
$$
\makdo=\frac{q^{-\nu^2+1/2}}{2\sqrt{a_\nu a_{-\nu}}\sqrt{z}}
E_q\left(-\frac{1-q^2}2z\right)\Phi_\nu(-z),
$$
and hence it is a holomorphic function in the region
$z\ne0$.
\end{predl}
\begin{predl}\label{p2.10}
The function $\makdj$ satisfies the relations
\beq{2.14}
\frac{2}{(1+q)z}\p_qz^\nu\makdj=-q^{-\frac{2-\de}2(\nu-1)}z^{\nu-1}
K_{\nu-1}^{(j)}((1-q^2)q^{\frac{2-\de}2}z;q^2),
\eq
\beq{2.15}
\frac{2}{(1+q)z}\p_qz^{-\nu}\makdj=-q^{\frac{2-\de}2(\nu+1)}z^{\nu-1}
K_{\nu+1}^{(j)}((1-q^2)q^{\frac{2-\de}2}z;q^2).
\eq
\end{predl}
\begin{predl}\label{p2.11}
The functions $\makdj$ satisfy the functional relations:
$$
q^{-\frac{2-\de}2\nu}K_{\nu-1}^{(j)}((1-q^2)z;q^2)-
q^{\frac{2-\de}2\nu}K_{\nu+1}^{(j)}((1-q^2)z;q^2)=
$$
\beq{2.16}
=-\frac{2}{(1-q^2)z}(q^{-\nu}-q^\nu)K_\nu^{(j)}((1-q^2)q^{\frac\de2}z;q^2),
\eq
$$
q^{-\frac{2-\de}2\nu}K_{\nu-1}^{(j)}((1-q^2)z;q^2)+
q^{\frac{2-\de}2\nu}K_{\nu+1}^{(j)}((1-q^2)z;q^2)=
$$
\beq{2.17}
-\frac{4}{(1-q^2)z}K_\nu^{(j)}((1-q^2)q^{\frac{\de-2}2}z;q^2),
+\frac{2}{(1-q^2)z}(q^{-\nu}+q^\nu)K_\nu^{(j)}((1-q^2)q^{\frac\de2}z;q^2).
\eq
\end{predl}
\begin{predl}\label{p2.12}
For any $\nu$ the functions $\mbesj$ and $\makdj$ form a fundamental
system of the solutions of equation (\ref{2.8})
\end{predl}

The $q$-Wronskian $W(I_\nu^{(j)},K_\nu^{(j)})(z)$ differs from the
$q$-Wronskian (\ref{2.9}) by constant multiplier.
$$
W(I_\nu^{(j)},K_\nu^{(j)})(z)=\left\{
\begin{array}{lcl}
\frac{q^{-\nu^2}(1-q^2)}2\sqrt{\frac{a_\nu}{a_{-\nu}}}
\qqexp & {\rm for}& j=1\\
\frac{q^{-\nu^2}(1-q^2)}2\sqrt{\frac{a_\nu}{a_{-\nu}}}
\qqeexp & {\rm for}& j=2.\\
\end{array}
\right.
$$
Hence, this function does not vanish.

\vspace{5mm}
{\bf 2.4. The $q$-integral representations}

\vspace{5mm}

In \cite{OR2} the representations of $q$-MFs (for $j=1, 2$) by the Jackson
$q$-integral were received. Let $z$ and $s$ be non commuting variables
$$
zs=qsz.
$$
\begin{predl}\label{p2.13}
If $\nu>\frac32$ the $q$-MFs can be represented by $q$-integrals
$$
\makd=\frac{q^{-\nu^2+\frac12}\G_{q^2}(\nu+\frac12)\G_{q^2}(\frac12)}
{4{\bf Q}_\nu}\sqrt{\frac{a_\nu}{a_{-\nu}}}
\left(\frac z2\right)^{-\nu}\times
$$
\beq{2.18}
\times\int_{-\infty}^\infty E_q\left(i\frac{1-q^2}2zs\right)
\frac{(-q^2s^2,q^2)_\infty}{(-q^{-2\nu+1}s^2,q^2)_\infty}d_qs,
\eq
$$
\makdo=\frac{q^{-\nu^2+\nu}\G_{q^2}(\nu+\frac12)\G_{q^2}(\frac12)}
{4{\bf Q}_{\frac12}}\sqrt{\frac{a_\nu}{a_{-\nu}}}
\left(\frac z2\right)^{-\nu}\times
$$
\beq{2.19}
\times\int_{-\infty}^\infty
\phantom._0\Phi_1\left(-;0;q,i\frac{1-q^2}2zs\right)
\frac{(-q^{2\nu+1}s^2,q^2)_\infty}{(-s^2,q^2)_\infty}d_qs,
\eq
where
\beq{2.20}
{\bf Q}_\nu=(1-q)\sum_{m=-\infty}^\infty\frac1{q^{m+\nu-\frac12}+
q^{-m-\nu+\frac12}}.
\eq
\end{predl}

\section{The $q$-Macdonald functions of kind 3}
\setcounter{equation}{0}

\begin{defi}\label{d3.1}
Starting from (\ref{2.18}), (\ref{2.19}) we define the $q$-Macdonald
function of kind 3 by its $q$-integral representation for $\nu>\frac32$.
$$
\makdon=C(\nu)\frac{q^{-\nu^2+\frac\nu2+\frac14}\G_{q^2}(\nu+\frac12)
\G_{q^2}(\frac12)}{4{\bf Q}_{\frac\nu2+\frac14}}
\left(\frac z2\right)^{-\nu}\times
$$
\beq{3.1}
\times\int_{-\infty}^\infty
\phantom._0\Phi_2\left(-;0,-q^{\frac12};q^{\frac12},-i\frac{1-q^2}2zs\right)
\frac{(-q^{\nu+\frac32}s^2,q^2)_\infty}
{(-q^{-\nu+\frac12}s^2,q^2)_\infty}d_qs.
\eq
Moreover
\beq{3.2}
\lim_{m\to\infty}\makdon|_{\frac{1-q^2}2z=q^m}=0.
\eq
\end{defi}

Multiplying both sides of (\ref{3.1}) by $(\frac z2)^\nu$ and setting
$z=0$, we obtain
$$
\left(\frac z2\right)^\nu\makdon|_{z=0}=C(\nu)
\frac{q^{-\nu^2+\frac\nu2+\frac14}\G_{q^2}(\nu+\frac12)\G^{q^2}(\frac12)}
{4{\bf Q}_{\frac\nu2+\frac14}}
\int_{-\infty}^\infty\frac{(-q^{\nu+\frac32}s^2,q^2)_\infty}
{(-q^{-\nu+\frac12}s^2,q^2)_\infty}d_qs.
$$
Calculate the last $q$-integral using the properties of the $q$-binomial
formula \cite{OR2}.
$$
Int=\int_{-\infty}^\infty\frac{(-q^{\nu+\frac32}s^2,q^2)_\infty}
{(-q^{-\nu+\frac12}s^2,q^2)_\infty}d_qs=2(1-q)\sum_{m=-\infty}^\infty q^m
\frac{(-q^{\nu+2m+\frac32},q^2)_\infty}{(-q^{-\nu+2m+\frac12},q^2)_\infty}=
$$
$$
=2(1-q)\sum_{m=-\infty}^\infty q^m\frac{(q^{2\nu+1},q^2)_\infty}
{(q^2,q^2)_\infty}\sum_{k=0}^\infty\frac{(q^{-2\nu+1},q^2)_kq^{(2\nu+1)k}}
{(q^2,q^2)_k(1+q^{-\nu+\frac12+2m+2k})}.
$$
The inner series converges uniformly with respect to $m$, and we can
change the order of summing. Then we have using (\ref{2.20})
$$
Int=\frac{2(1-q)(q^{2\nu+1},q^2)_\infty}{(q^2,q^2)_\infty}
\sum_{k=0}^\infty\frac{(q^{-2\nu+1},q^2)_kq^{2\nu k}}{(q^2,q^2)_k}
\sum_{m=-\infty}^\infty\frac{q^{m+k}}{1+q^{-\nu+\frac12+2m+2k})}=
$$
$$
=\frac{2q^{\frac\nu2-\frac14}(q^{2\nu+1},q^2)_\infty}{(q^2,q^2)_\infty}
{\bf Q}_{\frac\nu2+\frac14}
\sum_{k=0}^\infty\frac{(q^{-2\nu+1},q^2)_kq^{2\nu k}} {(q^2,q^2)_k}=
2q^{\frac\nu2-\frac14}{\bf Q}_{\frac\nu2+\frac14}
\frac{(q^{2\nu+1},q;q^2)_\infty}{(q^{2\nu},q^2;q^2)_\infty}.
$$
Hence
\beq{3.3}
\left(\frac z2\right)^\nu\makdon|_{z=0}=C(\nu)\G_{q^2}(\nu)
\frac{q^{-\nu^2+\nu}}2.
\eq

Let $\nu\ne n$. As $\makdon$ is the solution of (\ref{2.8}) for $\de=1$,
$\mbesse$ and $\mmbesse$ form the fundamental system of solutions of
this equation, and\\ $K_{-\nu}^{(3)}((1-q^2)z;q^2)=\makdon$ we can write
\beq{3.4}
\makdon=\frac12q^{-\nu^2+\nu}\G_{q^2}(\nu)\G_{q^2}(1-\nu)
[A_{-\nu}\mmbesse-A_\nu\mbesse].
\eq
It follows from this equality
\beq{3.5}
\left(\frac z2\right)^\nu\makdon|_{z=0}=
\frac12q^{-\nu^2+\nu}\G_{q^2}(\nu)A_{-\nu}.
\eq
It follows from (\ref{3.3}) and (\ref{3.5}) $C(\nu)=A_{-\nu}$.

Consider the restriction of function $\mbesse$ (\ref{2.7}) on the lattice
$\{q^n\}$. e.i. assume $\frac{1-q^2}2z=q^n$. Then
$$
\mbesse|_{\frac{1-q^2}2z=q^n}=I_\nu^{(3)}(2q^n;q^2)=\frac1{(q^2,q^2)_\infty}
\sum_{k=0}^\infty\frac{q^{k(\nu+k)+n(\nu+2k)}(q^{2\nu+2k+2},q^2)_\infty}
{(q^2,q^2)_k}=
$$
$$
=\frac{q^{\nu n}}{(q^2,q^2)_\infty}\sum_{k=0}^\infty
\frac{q^{k(k-1)+k(\nu+1+2n)}}{(q^2,q^2)_k}
\sum_{m=0}^\infty\frac{(-1)^mq^{m(m-1)}q^{(2\nu+2k+2)m}}{(q^2,q^2)_m}=
$$

$$
=\frac{q^{\nu n}}{(q^2,q^2)_\infty}
\sum_{m=0}^\infty\frac{(-1)^mq^{m(m+1)}q^{2\nu m}}{(q^2,q^2)_m}
\sum_{k=0}^\infty\frac{q^{k(k-1)+k(\nu+1+2m+2n)}}{(q^2,q^2)_k}=
$$
$$
=\frac{q^{\nu n}}{(q^2,q^2)_\infty}
\sum_{m=0}^\infty\frac{(-1)^mq^{m(m+1)}q^{2\nu m}}{(q^2,q^2)_m}
(-q^{\nu+1+2m+2n},q^2)_\infty=
$$
$$
=\frac{q^{\nu n}(-q^{\nu+1+2n},q^2)_\infty}{(q^2,q^2)_\infty}
\sum_{m=0}^\infty\frac{(-1)^mq^{m(m+1)}q^{2\nu m}}
{(q^2,q^2)_m(-q^{\nu+1+2n},q^2)_\infty}.
$$
Obviously the last sum tends to unit if $n\to-\infty$.

Now consider the limit of quotient
$$
\frac{\mbesse}{\mmbesse}\Bigr|_{\frac{1-q^2}2z=q^n}
$$
if $n\to-\infty$.
$$
\lim_{n\to-\infty}\frac{I_\nu^{(3)}(2q^n:q^2)}{I_{-\nu}^{(3)}(2q^n:q^2)}=
\lim_{n\to-\infty}\frac{q^{2\nu n}(-q^{\nu+1+2n},q^2)_\infty}
{(-q^{-\nu+1+2n},q^2)_\infty}=
\lim_{n\to\infty}\frac{(-q^{-\nu+1},q^2)_n(-q^{\nu+1},q^2)_\infty}
{(-q^{\nu+1},q^2)_n(-q^{-\nu+1},q^2)_\infty}=1.
$$
In order to (\ref{3.2}) is fulfilled it is necessary $A_{-\nu}=A_\nu$.
So $C(\nu)=A_{-\nu}=constant$. Assume this $constant$ equal to unit.
Then finally, if $\nu\ne n$
\beq{3.6}
\makdon=\frac12q^{-\nu^2+\nu}\G_{q^2}(\nu)\G_{q^2}(1-\nu)
[\mmbesse-\mbesse].
\eq
This definition must be extended to integral
values of $\nu=n$ by passing to the limit in (\ref{3.6}).

It is easily to verify a correctness of the following propositions.
\begin{predl}\label{p3.1}
The function $\makdon$ satisfies (\ref{2.14}), (\ref{2.15}).
\end{predl}
\begin{predl}\label{p3.2}
The function $\makdon$ satisfies functional relations (\ref{2.16}),
(\ref{2.17}).
\end{predl}
\begin{predl}\label{p3.3}
For any $\nu ~~\mbesse$ and $\makdon$ form the fundamental system of the
solutions of equation (\ref{2.8}) for $\de=1$, and
$$
W(I_\nu^{(3)},K_\nu^{(3)})=\frac12q^{-\nu^2}(1-q^2).
$$
\end{predl}
\begin{rem}\label{r3.1}
$$
\lim_{q\to1-0}\makdj=K_\nu(z).
$$
\end{rem}

It easily to show that $\lim_{q\to1}{\bf Q}_\nu=\frac\pi2$ for
(\ref{2.20}).
\begin{rem}\label{r3.2}
If $q\to1-0$ we have well-known integral representation for the
classical Macdonald functions
$$
K_\nu(z)=\frac{\G(\nu+\frac12)}{2\G(\frac12}\left(\frac z2\right)^{-\nu}
\int_{-\infty}^\infty(z^2+1)^{-\nu-\frac12}e^{izs}ds.
$$
\end{rem}

\small{
 }

\end{document}